\documentclass[a4paper, 11pt]{article} 
\usepackage{amsfonts,amsmath,amssymb,amscd}
\usepackage[all]{xy} 
\usepackage{latexsym}
\usepackage[T1]{fontenc}

\newtheorem{thm}{Theorem}[section]
\newtheorem{definition}{Definition}[section]
\newtheorem{prp}[thm]{Proposition}
\newtheorem{lem}[thm]{Lemma}
\newtheorem{cor}[thm]{Corollary}

\newenvironment{prf}[1][]{\begin{trivlist}\item[]{\bf Proof#1. }}
{\hfill $\blacksquare$ \end{trivlist}}

\newcommand{\nap}{\nabla^\perp}
\newcommand{\coker}{\mathrm{coker }}
\newcommand{\ind}{\mathrm{index}\,}

\newcommand{\bd}{\mathcal{D}}
\newcommand{\grad}{\vec{\nabla}}
\newcommand{\bfff}{\mathcal{F}}
\newcommand{\beee}{\mathcal{E}}

\newcommand{\br}{\mathcal{R}}
\newcommand{\bl}{\mathcal{L}}

\newcommand{\bm}{\mathcal{M}}

\newcommand{\ba}{\mathcal {A}}

\newcommand{\na}{\nabla}
\newcommand{\be}[1]{\begin{equation}\label{#1}}
\newcommand{\ee}{\end{equation}}

\newcommand{\kah}{\text{K\"{a}hler }}

\newcommand{\beq}{\begin{eqnarray*}}
\newcommand{\eeq}{\end{eqnarray*}}
\newcommand{\bpr}{\begin{prf}}
\newcommand{\epr}{\end{prf}}

\newcommand{\T}{\mathbb{T}}

\newcommand{\R}{\mathbb{R}}

\newcommand{\C}{\mathbb{C}}

\newcommand{\dsurdt}{\frac{\partial}{\partial t}}

\begin{document}
\title{Vanishing theorems for associative submanifolds}
\author{Damien Gayet}
\maketitle

\centerline{\textbf{Abstract}}
Let $M^7$ a manifold with holonomy in $G_2$, 
and $Y^3$ an associative submanifold with boundary in a coassociative
submanifold.  In \cite{GaWi}, the authors proved that $\bm_{X,Y}$, the moduli space of its associative deformations with boundary  in the fixed $X$, has finite virtual dimension. Using Bochner's technique, we give a vanishing theorem that forces
$\bm_{X,Y}$ to be locally smooth.

\bigskip

\textsc{MSC 2000:} 53C38 (35J55, 53C21, 58J32).

\smallskip

\textsc{Keywords:} $G_2$ holonomy; calibrated submanifolds; elliptic boundary problems ; Bochner technique
%
%
%
%
%
\section{Introduction}
Let $M$ a 7-dimensional riemannian manifold with holonomy included in $G_2$. In this case by parallel transport, $M$ supports a non degenerate harmonic  3-form $\phi$ with $\nabla \phi = 0$. Moreover,
 $M$ benefits a vector product $\times$ defined by 
$$<u\times v, w>= \phi(u,v,w),$$
even if to be spin is enough for the existence of this product. 
A 3-dimensional submanifold $Y$ is said
\textit{associative} if its tangent bundle is stable under the
vector product. In other terms, 
$\phi$ restricted to $Y$ is a volume form. 
Likely, a 4-dimensional submanifold $X$ 
is $coassociative$ if  the fibers of its
normal bundle are associative, or equivalently, 
$\phi_{|TX}$ vanishes. We refer to the abundant literature on this subject, see \cite{HaLa} or \cite{GaWi} for a summary with our notations.\\

\noindent
\textbf{The closed case.}
It is known from \cite{McL} that the deformation of
an associative submanifold $Y$ without boundary is an elliptic problem, and hence is of 
vanishing index. In general, the situation is obstructed. For instance,
consider  the flat torus $\T^3\times \{pt\}$ in the flat torus $\T^7 = \T^3\times \T^4$. 
This is an associative submanifold, and its moduli space $\bm_{\T^3\times \{pt\}}$ of associative
deformations contains at least the 4-dimensional $\T^4$.  \\

\noindent
A natural question is to find conditions which force the moduli space to be smooth, or
in other terms, which force the cokernel of the problem to vanish. For the closed case, 
Abkulut and Salur \cite{AkSa2} allow a certain freedom for the connection on the normal bundle,
the definition of associativity  
 and genericity. But examples
 are often non generic, and we would like to get a condition that is not a perturbative one. 
For holomorphic curves in dimension 4, there are topological conditions
on the degree of the normal bundle which imply smoothness
of the moduli space of complex deformations, see for example 
 \cite{HoLiSi}. The main point for this is that holomorphic
 curves intersect positively. In our case, there is no such phenomenon. \\

In \cite{McL}, page 30, McLean gives an example 
of an isolated associative submanifold. Since this was the start point of our work,
we recall it.  In \cite{BrSa},  Bryant and Salamon constructed  a metric of holonomy $G_2$ on the spin bundle $S^3\times \R^4$ of the round 3-sphere. In this case, the basis $S^3\times \{0\}$ is 
associative, the normal bundle is just the spin bundle of $S^3$, and the operator
related to the associative deformations of $S^3$ is just the Dirac operator. 
By the famous theorem of Lichnerowicz \cite{Li}, there are no non trivial harmonic spinors on $S^3$
for metric reasons (precisely, the riemaniann scalar curvature is positive), 
so the sphere is isolated as an associative submanifold. \\

\noindent 
\textbf{Minimal submanifolds.}
Recall that in manifolds with holonomy in $G_2$, 
associative submanifolds are minimal (the condition $d\phi=0$ is enough). 
In \cite{Si}, Simons gives a metric condition for a minimal submanifold
to be stable, i.e  isolated. For this, he introduce the following operator, a sort
of partial Ricci operator :
\begin{definition} Let $(M,g)$ a riemaniann manifold, $Y^p$ a submanifold in $M$ and $\nu$ its normal bundle. Choose $\{e_1,\cdots e_p\}$ 
a local orthonormal frame field of $TY$, and define the 0-order operator 
\beq
\br  : \Gamma(Y,\nu) &\longrightarrow& \Gamma(Y,\nu)\\
\psi &\mapsto & \pi_\nu \sum_{i=1}^p  R (e_i, \psi) e_i,
\eeq
where $R$ is the curvature tensor on $M$   and $\pi_\nu$
the orthogonal projection on $\nu$.
\end{definition}
\textit{Fact.} The definition is independant of the choosen oriented orthonormal frame, and 
$\br$ is symmetric.\\

\noindent
He introduces another operator $\mathcal A$ related to the second fondamental form of $Y$ :
\begin{definition} Let $SY$  the bundle over $Y$ whose fibre at a point $y$ is the space of symmetric endomorphisms
of $T_yY$, and $A \in Hom (\nu, SY)$ the second fundamental form defined by 
$$  A(\phi) (u) = -\nabla^\top_u \phi,$$ where $u\in TY$, $\phi \in \nu$, and $\nabla^\top$ is the projection on $TY$
of the  ambient Levi-Civita connection.
Consider the operator 
\beq
\mathcal A : \Gamma(Y,\nu) &\longrightarrow& \Gamma(Y,\nu)\\
\psi &\mapsto & A^t\circ A (\psi),
\eeq
where $A^t$ is the transpose of $A$.
\end{definition}
\textit{Fact.} This is a symmetric positive 0-th order operator. Moreover, it vanishes if $Y$ is totally geodesic. \\

\noindent
Using both operators, Simons gives a sufficient condition for a minimal submanifold 
to be stable :
\begin{thm}[\cite{Si}]\label{Simons}
Let $Y$ a minimal submanifold in $M$, and suppose that $\br - \mathcal A$ 
is positive. Then $Y$ cannot be deformed as a minimal submanifold. 
\end{thm}

\noindent
\textbf{Bochner technique.}
If $Y$ is an associative submanifold in $M$, 
we will recall that there is an operator $D$ acting on the normal vector fields of $Y$, 
such that its kernel can be identified with the infinitesimal associative deformations 
of $Y$. We will compute $D^2$ to use the Bochner method, and get vanishing theorems.
For this, we introduce the normal equivalent of the invariant second derivative. More
precisely, for every local vector fields $v$ and $w$ in $\Gamma(Y,TY)$, let 
 $$\nabla^{\perp 2}_{v,w} = \nap_v\nap_w - \nap_{\nabla^\top_v w},$$
acting on $\Gamma(Y,\nu)$. It is straightforward to see that it is tensorial in  $v$ and $w$.
 Moreover, define the equivalent of the connection laplacian : 
$$ \nabla^{\perp *}\nap = - \text{trace } (\nabla^{\perp 2}) = - \sum_i \nabla^{\perp 2}_{e_i,e_i}.$$ 
\begin{thm}\label{D^2}
For $Y$ an associative submanifold, 
$ D^2 =  \nabla^{\perp *}\nap + \br - \mathcal A.$
\end{thm}
\textit{Remark.} In fact this shows that for closed submanifolds, 
associativity does not give more conditions than the one for minimal submanifods, 
as long as we use this method. \\

\noindent
\textbf{The  case with boundary.}
In \cite{GaWi}, the authors proved that 
the deformation of an associative submanifold $Y$ with boundary in 
a coassociative submanifold $X$ is an elliptic problem of finite index. Moreover,
they  gave the value of this index in terms of a certain Cauchy-Riemann operator 
related to the complex geometry of the boundary. We sum up in the
following the principal results of the paper : 
\begin{thm}[\cite{GaWi}]\label{boundary}
Let $\nu_X$ the normal complementary of $T\partial Y$ in $TX_{|\partial Y}$, and 
$n$ the inner unit vector normal to $\partial Y$ in $Y$. Then 
the bundle $\nu_{X}$ is a subbundle of  $\nu_{|\partial Y}$ and is stable under the left action by $n$ under $\times$, as well as
the orthogonal complement $\mu_X$ of $\nu_X$  in $\nu$. Viewing  $T\partial Y$, $\nu_X$ and $\mu_X$ 
as $n\times$--complex line bundles, we have $\mu_X^*\cong\nu_X\otimes_{\C}T\partial Y$.
Besides, the problem of the associative deformations of $Y$ with boundary in $X$ 
is elliptic and of index
$$\ind (Y,X) = \ind \overline \partial_{\nu_X} = c_1(\nu_X) + 1- g,$$
where $g$ is the genus of $\partial Y$.
\end{thm}
\noindent
In this context, we introduce a new geometric object 
that is related to the geometry on the boundary :
\begin{prp}\label{delta}
Choose $\{v,w= n\times v\}$ a local orthonormal frame for $T\partial Y$.
Let $L$ a real plane subbundle of $\nu$ invariant by the action of $n\times$.
We define 
\beq
\bd_L : \Gamma(\partial Y,L) &\longrightarrow & \Gamma(\partial Y,L)\\
\phi &\mapsto &  \pi_L(v\times \nabla^\perp_w \phi- w\times \nabla^\perp_v \phi),
\eeq
where $\pi_L$ is the orthogonal projection on $L$ and $\nap$ the normal connection on $\nu$ induced
by the Levi-Civita connection $\nabla$ on $M$. Then 
$\bd_L$ is independant of the choosen oriented frame, 
of order 0  and symmetric.
\end{prp}
Now, we can express our main theorem :
\begin{thm}\label{open}
Let $Y$ an associative submanifold of a $G_2$-manifold $M$ with boundary in a coassociative $X$. 
If  $\bd_{\mu_X}$ and 
$\br-\mathcal A $ are positive, the moduli space $\bm_{Y,X}$ is locally smooth and of dimension 
given by the virtual one $\ind (Y,X)$.
\end{thm}
When $M = \R^7$, we get the following very explicit example considered in \cite{GaWi}. 
 Take a ball $Y$ in $\R^3\times \{0\}\subset \R^7$, with real analytic boundary, 
and choose $e$ any constant vector field in $\nu = Y \times \{0\}\oplus \R^4$.  By \cite{HaLa}, 
there is a unique local coassociative $X_e$ containing $\partial Y\times \R e$, such that 
$$TX_{|\partial Y} = T\partial Y \oplus \nu_X = T\partial Y \oplus \text{Vect }(e,n\times e).$$
Of course, the translation in the $e$-direction gives associative deformations of $Y$ with boundary in $X_e$. 
The next corollary shows that under a simple metric condition, this is the only way to deform $Y$ :
\begin{cor}\label{ball}
If $Y$ is a strictly convex ball in $\R^3$, then 
$\bm_{Y,X_e}= \R$. 
\end{cor}

\textbf{The Calabi-Yau extension. } 
Let $(N, J, \Omega, \omega)$ a Calabi-Yau 6-dimensional manifold, 
where $J$ is an integrable complex stucture, $\Omega$ a non vanishing holomorphic 3-form
and $\omega $ a \kah form. Here we allow holonomies which are only  subgroups of  $SU(3)$.
Then $M = N\times S^1$  
is a manifold with holonomy in $SU(3) \subset G_2$. The associated 
calibration 3-form is given by  $$ \phi = \omega\land dt + \Re \Omega. $$
Recall that a special lagrangian in $N$ is 3-dimensional submanifold $L$ in $N$ satisfying both 
conditions $\omega_{|TL} = 0$ and $\Im \Omega_{|TL} = 0$.  We know from \cite{McL} that $\bm_L$ the moduli space of  special lagrangian deformations of $L$ is smooth
and of dimension $b^1(L)$. Now every product   $Y= L\times \{pt\}$ 
of a special lagrangian and a point is an associative submanifold of $M$. \\

\noindent
If $\Sigma$ is a complex surface of $N$, then $X = \Sigma \times \{pt\}$
is a coassociative submanifold of $M$. Consider the problem of associative deformations
of  $Y = L\times \{pt\}$  with boundary in $X$ : 
\begin{thm}\label{SL-boundary} Let $L$ a special lagrangian submanifold in a 6-dimensional Calabi-Yau $N$,
such that $L$ has  boundary in a  complex surface $\Sigma$. Let $Y = L\times \{t_0\}$ in  $N\times S^1$
and $X = \Sigma \times \{t_0\}$. 
If the Ricci curvature of $L$ is positive, and the boundary of $L$ has positive mean curvature in $L$, 
 then $\bm_{Y, X}$ is locally smooth and has dimension $g$, where $g$ is the genus of $\partial L$. 
 \end{thm} 

 \textit{Aknowledgements. } I would like to thank the French Agence nationale pour la Recherche for its support, Vincent Borrelli (resp. Jean-Yves Welschinger) who convainced me that there is a life after curvature tensors (resp. Sobolev spaces), Gilles Carron and Alexei Kovalev for their interest in this work.
\section {Closed associative submanifolds}
\subsection{The operator $D$}
We begin with the version of McLean's theorem proposed by Akbulut and Salur. We will give below a new proof of this result.
\begin{thm}[\cite{McL},\cite{AkSa}]\label{Dirac} Let $Y$ an associative submanifold 
of a riemannian manifold $M$ with $G_2$-holonomy, and $\nu$ its normal bundle. Then the tangent space of its  associative deformations is the kernel of the operator 
\begin{eqnarray}\label{Operator}
D :  \Gamma(Y,\nu) &\longrightarrow & \Gamma(Y,\nu),  \nonumber \\
\psi &\mapsto &  \sum_{i=1}^3 e_i \times \nabla^\perp_{e_i} \psi,
\end{eqnarray}
where $\nabla^\perp$ is the connection on $\nu$ induced by the Levi-Civita connection $\nabla$ on $M$.
\end{thm}
As the proof of theorem \ref{D^2} is pretty technical, we refer to the last section for it.

\subsection{The implicit function machinery}
We give now the analytic elements needed for the clear definition of 
our moduli spaces and their smoothness.  
We define the Banach space $\beee$ of (not necessarily associative) embedding 
of a 3-manifold in $M$, and a function $F$ that measures the 
lack of associativity. Then we linearise $F$, and identify the tangent
space of the moduli space of associative deformations, a priori weak,  with its kernel. 
Since the derivative of $F$ is elliptic, those deformations are in fact smooth, and 
by the implicit function theorem, 
this moduli space is smooth if the cokernel vanishes.\\

Firstly, recall the existence on $(M,\phi)$ of a 
important object, the 3-form with values in $TM$ defined, for $u$, $v$, $w\in TM$
 by : 
\begin{eqnarray}\label{associator}
 \chi(u,v,w) = - u \times (v\times w) - <u,v>w + <v,w> v.
 \end{eqnarray}
  It is easy to check  \cite{AkSa} that  $\chi(u,v,w) $
 is orthogonal to the 3-plane $u\land v\land w$. Besides,
 if $Y$ is a 3-dimensional submanifold in $(M, \phi)$, 
 then $\chi_{|TY} = 0$ if and only if $Y$ is associative. \\

For the future computations, we will use the the following  usefull formula \cite{HaLa} :
$$ <\chi(u,v,w), \eta>= *\phi (u,v,w,\eta),$$
where  $*$ is the Hodge star, and $\eta \in TM$. So 
\begin{eqnarray}\label{star}
\chi  = \sum_k \eta_k \lrcorner *\phi \otimes \eta_k,
\end{eqnarray}
where $(\eta_k)_{k = 1, 2, \cdots 7} $ is
an orthonormal basis of the tangent space of $M$.\\

Now, as in \cite{McL}, we use this characterization 
to study the moduli space of associative deformations of
an associative $Y$.  

\bpr[ of theorem \ref{Dirac}] Let $(Y,g)$ any riemanian 3-manifod. For every embedding $ f : Y\to M$, define 
 $$F(f) = f^*\chi (\omega) \in  \Gamma (Y,  f_* TM),$$
where $\omega $ is the volume
form on $Y$. Then $f(Y)$ is associative if and only if  $F(f)$ vanishes. 
Consider a path of embeddings  $(f_t)_{t\in [0,1]}$. After  a reparametrization 
of $Y$, we can suppose that
$$s = \frac{d f_t}{d t}_{|t=0} \in \Gamma (Y,  \nu_{f_0}),$$
where $\nu_f$ is the normal bundle over $f(Y)$. \\

Suppose that $f_0(Y)$ is 
an associative submanifold of $M$, and that $f_0$ is
the injection of $Y$ in $ M$. 
In order to derivate the vector-valued form $F$ at $f_0$, we use the Levi-Civita connection :
$$ \nabla_{\dsurdt} F(f_t)_{|t=0} = \sum_k \bl_s (\eta_k\lrcorner *\phi) (\omega) \otimes \eta_k + (\eta_k\lrcorner *\phi) (\omega) \otimes \na_s \eta_k,$$
where $ \bl_s $ is the Lie derivative in the direction $s$.
The second member vanishes because $Y$ is associative.
Thanks to classical riemannian formulas, we compute the first term : 
 \beq
 \bl_s (\eta_k\lrcorner *\phi) (\omega)& =& (\eta_k\land \omega) \lrcorner \ \bl_s (*\phi) + ([\eta_k,s]\land \omega) \lrcorner *\phi \\
& = & \bl_s (*\phi) (\eta_k, \omega),
\eeq
since $ (\eta_k\lrcorner *\phi) (\omega) = <\chi(\omega) ,\eta_k> = 0.$
Writing $\omega = e_1\land e_2\land e_3$, where $(e_i)_{i = 1,2,3}$ 
is a local orthonormal frame of $TY$ with $e_3 = e_1\times e_2$, this is equal to 
$$
\nabla_s *\phi (\eta_k, \omega) + *\phi (\nabla_{\eta_k}s, \omega) + *\phi (\eta_k, \nabla_{e_1}s , e_2,e_3))+ *\phi (\eta_k, e_1,\nabla_{e_2}s ,e_3)) +
*\phi (\eta_k, e_1,e_2, \nabla_{e_3}s ).$$
The first term vanishes since $*\phi$ is covariantly constant, the second one vanishes too for the former reason,
and the third is 
$$*\phi (\eta_k, \nabla_{e_1}s , e_2,e_3) = *\phi (\eta_k, \nap_{e_1}s , e_2,e_3) = - <\nap_{e_1}s \times ( e_2 \times e_3),\eta_k>.$$
Using equation (\ref{star}) and the relation $e_2\times e_3 = e_1$ and summing up the two others similar terms,  
we have :
\begin{eqnarray}\label{diff}
Ds = \nabla_s F =  \sum_i e_i \times \nap_i s  .
\end{eqnarray}
\epr
After this linearisation, we come back to the problem of the moduli space. 
\begin{prp}\label{smoothness}
Let $Y$ an associative submanifold in $M$. If the kernel
of  the operator $D$ given by (\ref{Operator})  vanishes, then 
 $ \bm_Y $ is locally smooth and of vanishing dimension. 
If the moduli space $\bm_{Y}$ contains a smooth 
submanifold $\bm_k$ of dimension $k$, and if for every 
$Y'$ in $\bm_k$ the kernel of $D$ at $Y'$ is of dimension $k$,
then $ \bm_Y = \bm_k$. 
 \end{prp}
 \bpr 
For $kp>3$, it has a sense to consider the Banach space $$\beee = W^{k,p}(Y,M), $$
with tangent space at $f$ equal to $T_f \beee = W^{k,p}(Y,f_* TM)$.
Moreover for $(k-r)/3 > 1/p$, then 
$$W^{k,p}(Y,M) \subset C^{r} (Y,M), $$
and so $f\in \beee$ is $C^1$ if $k> 1+3/p$. 
In particular, one can define $\bfff$ the Banach bundle over $\bfff$ with fiber
$$ \mathcal{F}_{f} = W^{k-1,p}(Y,\nu_f).$$ 
It is clear that the operator $F$ extends to a section $F_{k,p}$ of $\bfff$ over $\beee$.
We just proved before that $F$ is differentiable, and if $f_0(Y) $ is associative, its derivative at $f_0(Y)$  along
a vector field $s \in T_{f_0}\beee$ is given by (\ref{diff}). \\

Now, the operator $D$ has symbol $$\sigma (\xi) : s \mapsto \sum_i \xi_i s \times e_i. = s\times \xi,$$
which is always inversible for $\xi \in TY\setminus \{0\}.$
This proves that $D$ is elliptic. 
Remark that $\sigma (\xi)^2s = - |\xi|^2 s$, which is the symbol
of the laplacian. 
Hence $\dim \ker D$ and $\dim \coker D$ are of finite dimension. 
By the implicit function theorem for Banach bundles, 
if $\coker D $ vanishes, then $F^{-1} (0)$ is a locally a  smooth Banach submanifold
of $\beee$ of finite dimension equal to  $\dim \ker D = \ind D$, 
which vanishes because $Y$ is odd-dimensional.
More generally, if $\dim \coker D$ is constant on 
the component of $\bm_X$ containing $Y$, then
the moduli space is still smooth of dimension $\dim \ker D$.
 Lastly, still because of ellipticity, all elements of $\bm_X$ are smooth.
\epr

\subsection{A vanishing theorem}
We can now formulate the following theorem, which can be deduced 
from theorem \ref{Simons}, since any associative submanifold is minimal.
\begin{thm}\label{closed}
Let $Y$ an associative submanifold of a $G_2$-manifold $M$. 
If the spectrum of $\br_\nu = \br - \mathcal A $ is positive, then 
$Y$ is isolated as an associative submanifold.
\end{thm}
 For reader's convenience, we give a proof of this result. 
 \bpr
Suppose that we are given a fixed closed associative submanifold $Y$. 
The virtual dimension of its moduli space of deformation
is vanishing. Consider a section $\psi \in \Gamma(Y,\nu)$. By classical calculations, using normal coordinates, 
we have 
\begin{eqnarray*}
-\frac{1}{2}\Delta |\psi|^2 & = & \sum_i<\nabla_i^\perp\psi, \nabla_i^\perp \psi> + <\psi, \nabla_i^\perp \nabla_i^\perp \psi>\\
                               & = & |\nabla^\perp \psi|^2 - <D^2\psi, \psi> + <\br_\nu \psi, \psi>
\end{eqnarray*}
by theorem \ref{D^2}.
 Since the laplacian is equal to $-\text{div}  (\grad)$, its integral over the closed $Y$ 
 vanishes. We get :
\begin{eqnarray}\label{integral}
0 = \int_Y |\nap \psi|^2 - <D^2 \psi, \psi > + <\br_\nu\psi, \psi> dy.
\end{eqnarray}
Suppose that  we have a section  $\psi\in \ker D$. Under the hypothesis that
 $\br_{\nu}$ is positive, the last equation implies  $\psi =0$. 
 Hence $ \dim \coker D = \dim \ker D = 0$, and by proposition \ref{smoothness}, 
 $\bm_Y$ is locally a smooth manifold of vanishing dimension, and $Y$ is isolated.  
\epr

\section{Associative submanifolds with boundary}
In this section we extend our result for rigidity in 
the case of associative submanifolds with boundary in 
a coassociative submanifold. In this case the index may be not zero, so 
rigidity transforms into smoothness of the moduli space.
\subsection{Implicit function machinery}
As before, define the adapted $\beee$, for $kp>3$ and $(k-r)/3 > 1/p$ : 
$$\beee_X = \{f \in W^{k,p} (Y,M), \ f(\partial Y)  \subset X \}.$$
This has the following tangent space :
$$T\beee_{X,f} = \{s \in W^{k,p} (Y,f^*TM), \  s_{|\partial Y}  \in f^*TX \}.$$ 
As before, we have the map :
$$ F : \beee_X \to W^{k-1,p} (Y,\nu_{f(Y)}).$$ 
It is enough to compute the derivative of $F_X$ at an application 
$f_0$ where $f_0(Y)$ is an associative submanifold. 
We suppose as in the closed case that $f_0$ is an injection $Y \hookrightarrow M.$ 
In this case,  lemma \ref{boundary} 
showed that $TX$ is orthogonal to $TY$ at $\partial Y$, hence
the derivative of $F$ at $f_0$ is :
$$ D : \{s \in W^{k,p} (Y,\nu),  \ s_{|\partial Y} \in \nu_X \} \to W^{k-1,p} (Y,\nu).$$ 
Now, to get some trace properties and use the results
of \cite{BoWo}, we need to restrict to the Sobolev
space $p = 2$. In particular, if $f\in H^s (Y,\nu) = W^{s,2}(Y,\nu)$, then $f_{|\partial Y} \in H^{s-\frac{1}{2}} (Y,\nu) $. 
By theorem 20.8 of \cite{BoWo}, 
the operator $D$ is Fredholm. In \cite{GaWi},
the authors computed its index, given by theorem \ref{boundary}.\\

\noindent
\textit{Notation} : For $L$ a subbundle of $\nu_{|\partial Y}$ of real rank equal to two,
define $$\ker (D, L) =  \{s \in W^{k,p} (Y,\nu),  \ s_{|\partial Y} \in L, \ ÊDs=0 \} .$$ 
 We will need the usefull
\begin{prp} \label{adjoint}
The operator $D$ is formally self-adjoint, i.e for $s$ and $s' \in \Gamma(Y,\nu)$,
\begin{equation}\label{self}
\int_Y < Ds,s'> - <s,Ds'> dy = - \int_{\partial Y}Ê< n\times s,s'> d\sigma, 
\end{equation}
where $d\sigma$ is the volume  induced by the
restriction of $g$ on the boundary, and $n$ is the normal inner unit vector of $\partial Y$. 
Moreover, $\coker (D,\nu_X) = \ker (D,\mu_X).$
\end{prp}
\bpr The proof of the firs assertion is \textit{mutatis mutandis}
the classical one for the classical Dirac operator, see
proposition 3.4 in \cite{BoWo} for example. For the reader's convenience
we give a proof of this. 
\beq 
 <Ds,s'> &= & <\sum_i e_i\times \nap_i s, s'> = -  \sum_i<\nap_i s,  e_i\times s'> \\
& = & -  \sum_i d_{e_i} < s,  e_i\times s'> + <s,\nap_i (e_i\times s')>\\
& = & -  \sum_i d_{e_i} < s,  e_i\times s'> + <s,\nabla_i^\top e_i\times s'+ e_i\times  \nap_i s'>.
\eeq
By a classical trick, define the vector field $X\in \Gamma (Y,TY)$ 
by $$ < X,w> = - <s, w\times s'> \ \forall w\in TY.$$
Note that the first product the one  of $TY$, 
and the second one the one of $\nu$. Now 
\beq - \sum_i d_{e_i} < s,  e_i\times s'>  &=& \sum_i d_{e_i} < X,  e_i> \\
& =&  \sum_i  <\nabla^\top_i X,  e_i> + < X,\nabla_i^\top e_i>\\
& =&  \sum_i   \text{div } X  -  <s,\nabla_i^\top e_i\times s'>.
\eeq 
By Stokes we get 
\beq
\int_Y<Ds,s'> dy &= & \int_{\partial Y} <X,-n> d\sigma + \int_Y <s,Ds'> dy \\
& =&  \int_{\partial Y} <s,n\times s'> d\sigma + \int_Y <s,Ds'> dy,
\eeq
which is what we wanted.
Now, let $s'\in \Gamma (Y,\nu)$ lying in $\coker (D,\nu_X)$. 
This is equivalent to say that for every $s\in \Gamma (Y,\nu_X)$,
we have $\displaystyle \int_Y <Ds,s'> dy= 0$. By the former result, 
we see that this equivalent to 
$$ \int_Y <s,Ds'> + \int_{\partial Y} <n\times s,s'> = 0.$$
This clearly implies that $Ds' = 0$, and $ s'_{|\partial Y} (\nu) \perp \nu_X,$
because $\nu_X$ is invariant under the action of $n\times$. 
So $s'_{|\partial Y} \in \mu_X$, and $s' \in \ker (D,\mu_X)$. 
The inverse inclusion holds too by similar reasons.
\epr

\subsection{Vanishing theorem}
\bpr[ of theorem \ref{open}] 
In order to get some smooth moduli spaces in the case with boundary,
we want to prove that $\coker (D,\nu_X)=  \ker (D,\mu_X)$ is trivial or has constant 
rank. So let $\psi \in   \ker (D,\mu_X).$
The boundary changes the integration (\ref{integral}), because the divergence has to be considered : 
\begin{eqnarray}
 \int_Y |\nap \psi|^2  + <\br_\nu \psi, \psi> dy = \frac{1}{2} \int_Y \text{div }\grad |\psi|^2 dy.
 \end{eqnarray}
By Stokes, the last is equal to 
\beq
                             -\frac{1}{2}\int_{\partial Y} d|\psi|^2(n) d\sigma =  -\int_{\partial Y} <\nap_n\psi, \psi> d\sigma,
                                                                  \eeq 
where $n$ is the normal inner unit vector of $\partial Y$. Choosing a local
orthonormal frame $\{v, w = n\times v\}$ of $T\partial Y$, and using the fact that $D\psi=0$, this is equal to  
\begin{eqnarray*}
                \int_{\partial Y} <w\times \nap_v \psi- v\times \nap_w\psi, \psi> d\sigma
           = -\int_{\partial Y} <\bd_{\mu_X} \psi,\psi> d\sigma.
\eeq
Summing up, we get the equation 
\begin{equation}\label{integration}
\int_Y |\nap \psi|^2 dy + \int_Y<\br_\nu \psi, \psi> dy + \int_{\partial Y} <\bd_{\mu_X} \psi,\psi> d\sigma=0.
\end{equation}
Now we can prove the theorem 1.4.
We see that if $\bd_{\mu_X}$ and $\br_\nu$ are
 positive, then $\psi$ vanishes. This means
that our deformation problem has no cokernel, 
and by a straigthforward generalization of proposition \ref{smoothness}, the moduli space is locally smooth.
\epr

\subsection{Some properties of the operator $\bd_L$}
We sum up the main results about $\bd_L$ in the following 
\begin{prp}\label{bd}
Let $Y$ an associative submanifold with boundary in a coassociative submanifold $X$, 
$L$ a subbundle of $\nu$ over $\partial Y$,  and $\bd_L$ as defined in the introduction. 
Then  $\bd_L$ is of order 0, symmetric, and its trace is $2H$, where
$H$ is the mean curvature of $\partial Y$ in $Y$ with respect to the outside
normal vector $-n$.
\end{prp}
\bpr Let  $L$ is a subbundle of $\nu$ invariant under the action of $n\times$.  
  It is straighforward to check that $\bd$ does not depend
of the orthonormal frame $\{v,w = n\times v\}$. For every $\psi\in \Gamma(\partial Y, L)$ and $f$ a function,
\beq
\bd_L(f\psi) &=& \pi_L (v\times \nabla_w (f\psi) - w\times \nabla_v (f\psi))\\
&= & f\bd_L \psi+ (d_w f)\pi_L ( v\times \psi) - (d_v f)\pi_L ( w\times \psi)  = f\bd_L \psi
\eeq
because $w\times L$ ans $v\times L$ are orthogonal to $L$. 
Now, decompose $$\nabla^\top = \nabla^{\top \partial} + \nabla^{\perp\partial}$$
into its two projections along  $T\partial Y$ and along the normal (in $TY$)  $n$-direction. 
For the computations, choose $v$ and $w = n\times v$ 
 the two orthogonal characteristic directions on $T\partial Y$, 
i.e  $\nabla_v^{\top\partial} n = - k_v v$
 and $\nabla_w^{\top\partial}  n =- k_w w$, where $k_v$ and $k_w$ are the two principal curvatures. 
 We have $\nabla^{\perp \partial}_v  =  k_v n$ and
 $<\nabla^{\perp\partial} _w v,n> = 0$, and the same, \textit{mutatis mutandis}, for $w$. 
Then, for $\psi$ and  $\phi\in  \Gamma(\partial Y, L)$,  using the fact that $T\partial Y\times L $
is orthogonal to $ L$, 
\beq 
<\bd_L \psi, \phi> &=& < \nabla^\perp_w (v\times \psi) - (\nabla_w^{\perp\partial} v)\times \psi- 
\nabla_v^\perp (w\times \psi) + (\nabla^{\perp\partial}_v w)\times \psi,\phi> \\
&=& < \nabla^\perp_w (v\times \psi) - \nabla^\perp_v (w\times \psi) ,\phi> \\
&=& - <v\times \psi, \nabla^\perp_w \phi> +  <w\times \psi ,\nabla^\perp_v\phi> \\
&= & <\psi, v\times \nabla^\perp_w \phi - w\times \nabla^\perp_v \phi> =  <\psi, \bd_L \phi>.
\eeq
To prove that the trace of $\bd_L$ is $2H$, let $e\in L$ a local unit section of $L$. 
We have $n\times e \in L$ too, and 
\beq
<\bd_L(n\times e),n\times e> &=& < v\times ((\nabla^{\top\partial}_w n)\times e) + v \times (n\times \nap_w e), n\times e> \\
&& - <w\times (\nabla^{\top \partial}_v n)\times e - w \times (n\times \nabla^\perp_v e), n\times e> \\
&=& <v\times (- k_w w \times e) -  w\times (-k_v v\times  e) , n\times e>\\
&&+ <v \times (n\times \nabla^\perp_w e) - w \times (n\times \nabla^\perp_v e), n\times e>\\
&=&  k_w + k_v - <n\times (w \times (n\times \nabla^\perp_v e) - v \times (n\times \nabla^\perp_w e)), e>\\
&=& 2H - <\bd_L e, e>.
\eeq
This shows that $\text{trace } \bd_L = 2H$. 
\epr

\subsection{Flatland}
In flat spaces, $R$ vanishes, and so  $\br_\nu = -\ba$. Hence a priori theorem \ref{open} does'nt apply. Nevertheless, we have the 
\begin{cor}
Let $Y$ a totally geodesic associative submanifold in a flat $M$, with boundary
in a coassociative $X$. If  $\bd_{\mu_X}$ 
 positive, then $\bm_{Y,X}$ is locally smooth and of expected dimension.
\end{cor}
\bpr The hypothesis on $Y$ implies that $\br_\nu =0$. 
Formula (\ref{integration}) shows 
that $\nabla^\perp \psi = 0$ and $\psi_{|\partial Y} = 0$. 
Using $d|\psi|^2 = 2<\nap \psi , \psi> = 0$, we get that $\psi = 0$, and $\coker (Y,\nu_X) = \ker (Y,\mu_X)=0$. 
\epr
\bpr[ of corollary \ref{ball}] Let $Y$ in $\R^3\times\{ 0\} \subset \R^7$, and $e\in \{0\}\times \R^4$.
From \cite{HaLa} the boundary of $Y$ lies in a local coassociatif
submanifold $X_e$ of $\R^7$, which contains $\partial Y\times \R e$ and whose
 tangent space over $\partial Y$ is
$T\partial Y\oplus \R e\oplus \R n\times e$.
We see that
$Y$ has a direction of associative deformation along the fixed direction $e$, hence
 the dimension of the kernel of our problem is bigger than 1. On the other hand, 
 the index is $c_1(\nu_X)+1-g = 1$.
 We want to show that $\bd_{\mu_X}$ is positive. To see that, 
 we choose orthogonal characteristic directions $v$ and $w = n\times v$ 
 in $T\partial Y$ as before.
 From theorem \ref{boundary}, we know that  $v\times e$ is a non vanishing section of $\mu_X$. 
 We compute :
\beq
\bd_{\mu_X}(v\times e) &  = & v \times ( \nabla^{\perp \partial}_w v \times e) - 
                                                     w  \times (\nabla^{\perp \partial}_v v \times e)\\
                                         &  = & - k_v w\times (n\times e) =  k_v v \times e.
\eeq
This shows that $k_v$ is an eigenvalue of $\bd_{\mu_X}$, and 
since we know that its trace is $2H$, we get that the other eigenvalue is $k_w$. 
Those eigenvalues are positive if the boundary of $Y$ is strictly convex. By the last corollary,
we get the result.
\epr
\textit{Remark. } In fact, we can give a better statement.  Indeed,
let $\psi \in  \ker (D,\nu_X)$, 
and decompose  $\psi_{|\partial Y}$ as $\psi = \psi_1 e + \psi_2 n\times e$.  
Of course, $e$ is in the kernel of $\bd_{\nu_X}$, and hence 
by proposition \ref{delta}, the second term is an eigenvector of $\bd_{\nu_X}$
for the eigenvalue $2H$. So
 formula (\ref{integration}) gives
$$ \int_{Y} |\nap \psi|^2+ \int_{\partial Y} 2H|\psi_2|^2 = 0.$$
If $H>0$, this imply immediatly that $\psi_2=0$, and $\psi_1$
is constant, so $\psi$ is proportional to $e$. This proves that
$\dim \ker (D,\nu_X) = 1$ under the weaker condition that $H>0$. 
 
\section{Extensions from the Calabi-Yau world}
\textbf{Closed extension}.
Let $(N, J, \Omega, \omega)$ a 6-dimensional manifold with holonomy in $SU(3)$. Then 
$M = N\times S^1$ is a manifold with holonomy in $G_2$, with 
the calibration the 3-form given by  $ \phi = \omega\land dt + \Re \Omega. $
Let $L$ a special lagrangian 3-dimensional submanifold in $N$. Recall that since $L$ is lagrangian, 
its normal bundle is simply $JTL$. 
Then $Y =L\times \{pt\}$ is an associative submanifold of $N\times S^1$, and
 its normal bundle $\nu$ is isomorphic to $JTL\times \R \partial_t$,
 where $\partial_t$ is the dual vector field of $dt$. 
Since the translation along $S^1$ preserves the associativity of $Y$, we hence
have $\bm_L\times S^1 \subset \bm_Y$.  We prove that
in fact,  there is equality, without any condition (compare an equivalent 
result for coassociative submanifolds by Leung in $\cite{Le}$) :
\begin{thm}\label{Calabi}
The moduli space $\bm_{L\times \{pt\}}$of associative deformations of $L\times \{pt\} $ is
always smooth, and can be identified with
the product $\bm_L\times S^1$. 
  \end{thm}

\bpr  In this situation, we don't use the former expression of $D^2$.
Instead, we give another formula for it. 
 If $s= J\sigma\oplus \tau \partial_t$ is a section of $\nu$, with $\sigma \in \Gamma(L,TL)$
 and $\tau \in \Gamma(L,\R)= \Omega^0(L)$, 
 we call $\sigma^\vee \in \Omega^1(L,\R)$ the 1-form dual to $\sigma$, and we use
 the same symbol for its inverse. Moreover, we use the classical notation  $* : \Omega^k(L) \to \Omega^{3-k}(L)$
 for  the Hodge star. Lastly, we define : 
 \beq 
 D^\vee : \Omega^1(L)\times \Omega^0(L) &\longrightarrow &\Omega^1(L)\times \Omega^0(L)\\
 (\alpha, \tau)& \mapsto& ((-J\pi_L D(J\alpha^\vee, \tau))^\vee, \pi_t D(J\alpha^\vee, \tau)),
 \eeq
  where $\pi_L$ (resp. $\pi_t$) 
 is the orthogonal projection $\nu = NL \oplus \R$ on the first (resp. the second) component. This is just a
 way to use  forms on $L$ instead of normal ambient vector fields.
\begin{prp}\label{harmonic}
For every $(\alpha, \tau) \in \Omega^1(L)\times \Omega^0(L)$, 
\beq 
D^\vee (\alpha,\tau) &=& (-*d\alpha - d\tau, *d*\alpha) \ \text{and} \\
D^{\vee 2} (\alpha, \tau) &=& - \Delta (\alpha, \tau),
\eeq
where $\Delta  = d^*d + dd^*$ (note that it is $d^*d$ on $\tau$). 
\end{prp}
Assuming for a while this propositioin, we see that for an infinitesimal associative  deformation of 
$L\times \{pt\}$, then $\alpha$ and $\tau$ are harmonic 
over the compact $L$.  In particular,  $\tau$ is constant and $\alpha$
describes an infinitesimal  special lagrangian deformation of $L$ (see \cite{McL}).  
In other words, the only way to displace $Y$ is 
to perturb $L$ as special Lagrangian  in $N$ and translate
it along the $S^1$-direction. Lastly, $\dim \coker D = \dim \ker D= b^1(L) + 1 $ and by proposition \ref{smoothness},
$ \bm_Y$ is smooth and of dimension $b^1(L) + 1$.  \epr

\begin{prf}[ of proposition \ref{harmonic}] We will use the simple  formula $ \nabla^\perp Js = J\nabla^\top s$ for all sections $s\in \Gamma(L, NL)$.
For $(s,\tau) \in \Gamma(L, NL)\times \R$, and $e_i$ local orthonormal frame on $L$, 
\beq 
D(s,\tau) &=& \sum_{i,j} <e_i\times \nabla^\perp_i s,Je_j>Je_j +   
 \sum_i<e_i\times \nabla^\perp_i s, \partial_t>\partial_t + 
 \sum_i \partial_i \tau \ e_i \times \partial_t\\
     &=& J\sum_{i,j} \phi(e_i,\nabla^\perp_i s, Je_j) e_j 
           + \sum_i \phi(e_i, \nabla^\perp_i s, \partial_t)\partial_t+
           J\sum_{i,j} \partial_i \tau \ <e_i\times \partial_t,Je_j>e_j,
           \eeq
           where we used that $e_i\times \partial_t \perp \partial_t$. 
           \beq
     &=& J\sum_{i,j}\Re \Omega(e_i, \nabla^\perp_i s,Je_j) e_j + 
                  \sum_i \omega(e_i,\nabla^\perp_i s) \partial_t+ 
                J  \sum_{i,j} \partial_i \tau \ \phi(e_i,\partial_t,Je_j)e_j \\
     &= & J\sum_{i,j} \Re \Omega (e_i, J\nabla_i^\top\sigma, Je_j)e_j +
                                \sum_i \omega(e_i, J\nabla_i^\top \sigma) \partial_t +
                          J      \sum_{i,j} \partial_i \tau \ \omega (Je_j, e_i) e_j,
                                \eeq
 where $\sigma = -Js \in \Gamma(L,TL)$.
                                \beq
     &= & -J \sum_{i,j} \Re \Omega (e_i, \nabla_i^\top\sigma, e_j)e_j + 
                                \sum_i <e_i,\nabla_i^\top \sigma> \partial_t -
                        J   \sum_{i,j} \partial_i \tau <e_j, e_i> e_j\\                          
     &= & -J \sum_{i,j} Vol (e_i, \nabla_i^\top\sigma, e_j) e_j + 
                                \sum_i <e_i,\nabla_i^\top \sigma> \partial_t - 
                                                          J \sum_{i} \partial_i \tau e_i  ,                            
                                                             \eeq
    since $\Re \Omega $ is the volume form on $TL$.
It is easy to find that this is equivalent to  
$$D(s,\tau) =  -J (*d\sigma^\vee)^\vee + (*d*\sigma^\vee) \partial_t - J (d\tau)^\vee,$$
and so $$D^\vee (\sigma^\vee, \tau) = (-*d \sigma^\vee - d\tau, *d*\sigma^\vee).$$
Now, since $d^* = (-1)^{3p+1}*d*$ on the $p$-forms, one easy checks the formula for $D^2$. 
 \end{prf}
\bpr[ of theorem \ref{SL-boundary}] 
Consider $L$ a special lagrangian with boundary in a complex surface $\Sigma$, 
and $Y = L\times \{pt\}$ (resp. $X = \Sigma\times \{pt\}$ )
its associative (resp. coassociative) extension.
It is clear that $\nu_X$ is equal (as a real bundle) to $JT\partial L$, and 
$\mu_X$ it the trivial $n\times$-bundle generated by $\partial_t$. 
We begin by computing the index of the boundary problem. 
This is very easy, since $\mu_X$ is trivial, and by theorem \ref{boundary}, 
we have $\nu_X \cong T\partial L^*$ (as $n\times$-bundles. Hence the index is 
$$-c_1(T\partial L) + 1-g = - (2-g)+1- g = g-1,$$
where $g$ is the genus of $\partial Y$.
 Now let $\psi =s+ \tau  \frac{\partial }{\partial t}  $
belonging to $\coker (D,\nu_X) = \ker (D, \mu_X)$, where
 $s$ a section of $NL$ and $\tau \in \Gamma(L,\R)$. 
Let $\alpha = -Js^\vee$. By proposition \ref{harmonic}, 
$\alpha$ is a harmonic 1-form, and $\tau$ is harmonic (note that 
$Y$ is note closed, so $\tau$ may be not constant). By classical results for harmonic 1-forms, we have : 
$$\frac{1}{2}\Delta |\psi|^2 =  \frac{1}{2}\Delta (|\alpha|^2 + |\tau|^2)= |\nabla_L \alpha|^2 + |d\tau|^2+ \frac{1}{2}\text{Ric } (\alpha, \alpha). $$
Integrating on $L\times \{pt\}$, we obtain the equivalence of formula (\ref{integration}) :
$$ - \int_{\partial Y} <\bd_{\mu_X} \psi,\psi> d\sigma=  \int_Y  |\nabla_L \alpha|^2 + |d\tau|^2+ \frac{1}{2}\text{Ric } (\alpha, \alpha) dy.$$
Lastly, let us compute the eigenvalues of $\bd_{\mu_X}. $
The constant vector $\frac{\partial }{\partial t}$ over $\partial Y$
lies clearly in the kernel of $\bd_{\mu_X}$. 
By proposition  \ref{bd}, the other eigenvalue of $\bd_{\mu_X}$ 
is $2H$, with eigenspace generated by $n\times \frac{\partial }{\partial t}$. Over $\partial Y$,  $s$ lies in $JTL\cap \mu_X$, hence is proportional
to $n\times \frac{\partial }{\partial t}$. Consequently, 
$\bd_{\mu_X} \psi  =  2 H s$ and  
$$ - \int_{\partial Y} 2H|s|^2 d\sigma=  \int_Y  |\nabla_L \alpha|^2 + |d\tau|^2+ \frac{1}{2}\text{Ric } (\alpha, \alpha) dy.$$

This equation, the  positivity of the Ricci curvature and the positivity of  $H$ show  that $\alpha$ vanishes and $\tau$ is constant. 
So we see that $\dim \coker (Y,X) = 1$, and by the constant rank theorem, $\bm_{Y,X}$ is locally smooth and of dimension
$\dim \ker (Y,X) = g$. 
\epr
\section{Computation of $D^2$}
\bpr[ of theorem \ref{D^2}]
Before diving into the calculi, we need the following trivial lemma :
\begin{lem}\label{lemme}
Let $\nabla$  the Levi-Civita connection on $M$ and $R$ its curvature tensor. 
For any vector fields $w$, $z$, $u$ and $v$ on $M$,
we have 
\beq 
\nabla (u\times v) &=& \nabla u \times  v + u \times \nabla v \text{ and }\\
R(w,z) (u\times v) &=& R(w,z)u\times v + u\times R(w,z)v.
\eeq
If $Y$ is an associative submanifold of $M$ with normal bundle $\nu$,  $u\in \Gamma (Y, TY)$, $v\in \Gamma (Y, TY)$ and $\eta \in \Gamma(Y, \nu)$,
then 
\beq
 \nabla^\top (u\times v) &=& \nabla^\top u \times  v + u \times \nabla^\top v \text{ and }\\
 \nabla^\perp (u\times \eta) &=& \nabla^\top u \times  v + u \times \nabla^\perp v,
\eeq
where $\nabla^\top= \nabla - \nap$ is the orthogonal projection of $\nabla$ on $TY$. 
\end{lem}
\bpr
Let $x_1, \cdots, x_7$ normal coordinates on $M$ 
near $x$, and $e_i = \frac{\partial}{\partial x_{i}}$
their derivatives,  orthonormal at $x$.
We have $$u\times v = \sum_i <u\times v, e_i> e_i = \sum_i \phi(u,v,e_i)e_i,$$
so that at $x$, where $\nabla_{e_j} e_i = 0$,
\begin{eqnarray*}
\nabla (u\times v) &= &\sum_i (\nabla \phi (u,v,e_i) + 
\phi (\nabla u,v,e_i) +  \phi ( u,\nabla v,e_i) 
+ \phi (\nabla u,v, \nabla e_i) )e_i\\
&=& \sum_i ( \phi (\nabla u,v,e_i)+ \phi ( u,\nabla v,e_i)) e_i =  
\nabla u \times  v + u \times \nabla v,
\end{eqnarray*}
because $\nabla \phi = 0$. Now if $u$ and $v$ are in $TY$, 
then we get the result after remarking that 
$(\nabla u \times v )^\top =   \nabla^\top u \times v$,
because $TY$ is invariant under $\times$. The last relation is implied by $TY\times \nu \subset \nu$
and $\nu\times \nu \subset TY$. The curvature relation 
is easily derived from the definition $R(w,z) = \nabla_w\nabla_z - \nabla_z\nabla_w - \nabla_{[w,z]}$
and the derivation of the vector product. 
\epr 
We compute $D^2$ at a point $x\in Y$. For this, we choose normal coordinates on $Y$ and $e_i\in \Gamma(Y,TY)$ their associated 
derivatives, orthonormal at $x$. To be explicit, $\nabla^\top e_i = 0$ at $x$.
Let $\psi \in \Gamma (Y,\nu). $     
      \begin{eqnarray*}
          D^2\psi &= & \sum_{i,j} e_i\times \nabla^\perp_i (e_j \times  \nabla^\perp_j \psi) \\  
         & = & \sum_{i,j} e_i\times (e_j \times \nabla^\perp_i \nabla^\perp_j \psi) 
          +\sum_{i,j} e_i\times (\nabla_i^\top e_j \times  \nabla^\perp_j\psi) 
          \\
          &=&
          -\sum_i\nabla^\perp_i \nabla^\perp_i \psi - \sum_{i\not= j} 
(e_i \times  e_j) \times \nabla^\perp_{i}\nabla^\perp_{j} \psi\\
          &=&
   \nabla^{\perp *}\nap \psi - 
  \sum_{i< j} (e_i \times  e_j) \times 
(\nabla^\perp_{i}\nabla^\perp_{j} -  \nabla^\perp_{j}\nabla^\perp_{i}) \psi\\
  &=&   \nabla^{\perp *}\nap \psi - 
  \sum_{i< j} (e_i \times  e_j) \times R^\perp(e_i,e_j) \psi.
  \eeq
  Since $(e_i \times  e_j) \times R^\perp(e_i,e_j)$ is symmetric in $i, j$, this is equal to 
 $$  \nabla^{\perp *}\nap \psi - 
  \frac{1}{2}\sum_{i, j} (e_i \times  e_j) \times R^\perp(e_i,e_j) \psi.
$$
The main tool for the sequence is the Ricci equation. Let $u$, $v$ in $\Gamma(Y,TY)$
and $\phi$, $\psi$   in   $ \Gamma(Y,\nu)$.  
$$ <R^\perp(u,v)\psi,\phi> =  <R(u,v)\psi,\phi> + < (A_\psi A_\phi- A_\phi A_\psi) u,v>,$$
where $A_\phi (u) = A(\phi)(u)  = -\nabla^\top_u \phi.$
Choosing $\eta_1, \cdots, \eta_4$ an
orthonormal basis of $\nu$ at the point $x$, we get 
\beq
-  \frac{1}{2}\sum_{i, j} (e_i \times  e_j) \times R^\perp(e_i,e_j) \psi &=& - \frac{1}{2}\sum_{i, j, k} <(e_i \times  e_j) \times R^\perp(e_i,e_j) \psi,\eta_k>\eta_k\\
&=& \frac{1}{2} \sum_{i, j, k} <R^\perp(e_i,e_j) \psi,(e_i \times  e_j) \times \eta_k>\eta_k\\
&=& -\frac{1}{2} \pi_\nu \sum_{i, j} (e_i \times  e_j) \times R(e_i,e_j) \psi \\
&&+
\frac{1}{2}\sum_{i, j, k}<(A_\psi A_{(e_i \times  e_j)\times\eta_k}- A_{(e_i \times  e_j)\times\eta_k} A_\psi) e_i , e_j>\eta_k.
\eeq
Using the classical Bianchi relation $R(e_i,e_j)\psi = -R(\psi, e_i)e_j - R(e_j,\psi) e_i$,
the first part of the sum is equal to 
\beq
I = - 2\pi_\nu (e_1\times R(e_2,\psi)e_3 + e_2\times R(e_3,\psi) e_1 + e_3 \times R(e_1,\psi)e_2) =\\
- 2\pi_\nu (e_1\times R(e_2,\psi)(e_1\times e_2) + e_2\times R(e_3,\psi) (e_2\times e_3) + 
e_3 \times R(e_1,\psi)(e_3\times e_1) )=\\
- 2\pi_\nu (e_1\times (R(e_2,\psi)e_1\times e_2+ e_1\times R(e_2,\psi)e_2) + 
                e_2\times (R(e_3,\psi) e_2\times e_3+  e_2\times R(e_3,\psi)e_1)+\\
                e_3\times (R(e_1,\psi) e_3\times e_1+  e_3\times R(e_1,\psi)e_2)) =\\
                - I + 2\pi_\nu \sum_i R(e_i,\psi) e_i,
\eeq
which gives $I = \pi_\nu \sum_i R(e_i,\psi) e_i$.\\

\noindent
The Weingarten endomorphisms are symmetric, so that the second part of the  sum is 
\beq 
\frac{1}{2}\sum_{i, j,k}< A_{(e_i \times  e_j)\times\eta_k}e_i,A_\psi e_j> \eta_k
- 
\frac{1}{2}\sum_{i, j,k}< A_\psi e_i,A_{(e_i \times  e_j)\times\eta_k}e_j> \eta_k.
\eeq
It is easy to see that the second sum is the opposite of the first one. 
We compute 
\beq 
A_{(e_i \times  e_j)\times\eta_k}e_i &=&-(\nabla^\perp_i e_i \times e_j)\times \eta_k
 - (e_i \times \nabla^\perp_i e_j)\times \eta_k
  + (e_i\times e_j)\times A_{\eta_k} e_i.
\eeq
But we know that an associative submanifold is minimal, so that
$$ \sum_i \nabla^\perp_i e_i =0.$$
Moreover, deriving the relation $e_3 = \pm e_1\times e_2$, one 
easily check that $$ \sum_i e_i\times \nabla^\perp_j e_i=0.$$
Summing, the only resting term is 
$$
\sum_{i, j, k}< (e_i\times e_j)\times A_{\eta_k}e_i,A_\psi e_j> \eta_k. $$
We now use the classical formula for vectors $u$, $v$ and $w$ in $TY$ :  
$$ (v\times w) \times u = <u,v>w - <u,w> v, $$
hence $$(e_i\times e_j)\times A_{\eta_k}e_i = <A_{\eta_k}e_i , e_i> e_j - <A_{\eta_k}e_i , e_j> e_i.$$
One more simplification comes from  $\sum_i < A_{\eta_k}e_i,e_i> =0$ 
for all $k$ because since $Y$ is minimal, 
so our sum is now equal to 
\beq
-\sum_{i, j, k}<A_{\eta_k}e_i , e_j> <e_i,  A_\psi e_j> \eta_k = -\mathcal{A}\psi.
\eeq
\epr


\bigskip

\noindent 
D.~\textsc{Gayet}\\
Universit\'e de Lyon, CNRS, Universit\'e Lyon 1, Institut Camille Jordan,\\
 F--69622 Villeurbanne Cedex, France\\
e-mail: \texttt{gayet@math.univ-lyon1.fr}

\medskip

\end{document}